\documentclass[a4paper, 12pt]{article}
\usepackage{amsfonts}
\usepackage{fancyhdr}
\usepackage{graphicx}
\usepackage{epsfig}
\usepackage{amssymb}
 \usepackage{mathrsfs,amsfonts,amsmath}
 \usepackage{color}

\makeatletter

\@addtoreset{equation}{section} \makeatother
\newtheorem{Theorem}{Theorem}[section]
\newtheorem{Remark}[Theorem]{Remark}
\newtheorem{Lemma}[Theorem]{Lemma}

\begin{document}
\title{\textbf{Strong convergence rates of modified truncated EM method for stochastic differential equations
\footnote{Supported by Natural Science Foundation of China (NSFC
11601025).}}}
\author{ Guangqiang Lan\footnote{Corresponding author: Email:
langq@mail.buct.edu.cn.}\quad and\quad
 Fang Xia
\\ \small School of Science, Beijing University of Chemical Technology, Beijing 100029, China}

\date{}

\maketitle

\begin{abstract}
Motivated by truncated EM method introduced by Mao (2015), a new explicit numerical method named modified truncated Euler-Maruyama method is developed in this paper. Strong convergence rates of the given numerical scheme to the exact solutions to stochastic differential equations are investigated under given conditions in this paper. Compared with truncated EM method, the given numerical simulation strongly converges to the exact solution at fixed time $T$ and over a time interval $[0,T]$ under weaker sufficient conditions. Meanwhile, the convergence rates are also obtained for both cases. Two examples are provided to support our conclusions.
\end{abstract}



\section{Introduction}

\noindent
Numerical methods for stochastic differential equations (SDEs) have been playing more and more important roles because most equations can not be solved explicitly. In general, there are two kinds of numerical methods, the one is explicit and the other is implicit. The most commonly used explicit numerical method is the well known Euler-Maruyama (EM) method. There are a lot of literature concerning with this method, e.g., \cite{Maruyama,Milstein,Kloeden,mao2,Hairer}. However, as mentioned in \cite{mao}, most of the existing strong convergence theory for numerical methods requires the coefficients of the SDEs to be globally Lipschitz continuous(see e.g. \cite{Kloeden,mao2}). In 2002, Higham et al. \cite{HMS} studied the strong convergence for numerical approximations under local Lipschitz condition for the first time plus the bounded condition on the $p$th moments of both exact and numerical solutions to the underlying SDE. Recently, Hutzenthaler et al. \cite{hut} proved, for a large class of SDEs with superlinearly growing coefficient functions, that both the distance in the strong $L^p$-sense and the distance between the $p$th absolute moments of the Euler approximation and of the exact solution of the SDE diverge to infinity for all $p\ge1$. Therefore, implicit methods have naturally been used to study the solutions to SDEs without the linear growth condition.

Implicit methods, including backward EM scheme, split-step backward scheme and $\theta$-EM scheme have been extensively studied. For example, Higham et al. \cite{HMS} studied convergence of a split-step backward Euler method for nonlinear SDEs under the assumption that the drift  satisfies the one-sided Lipschitz condition and the diffusion is globally Lipschitz, Mao and Szpruch \cite{MS} studied strong convergence rates for backward EM scheme for non-linear dissipative-type SDEs with super-linear diffusion coefficients, \cite{BT} analyzed strong convergence of split-step backward Euler method for SDEs with non-smooth drift.

Although \cite{hut} showed strong and weak divergence in finite time of Euler¡¯s method for stochastic differential equations under non-globally Lipschitz condition, some modified EM methods have recently been developed since they have simple algebraic structures, cheap computational costs and acceptable convergence rates. Recently, in \cite{mao}, Mao developed a new explicit numerical simulation method, called truncated EM method. Strong convergence theory were established there under local Lipschitz condition plus the Khasminskii-type condition. And then he obtained sufficient conditions for the strong convergence rate of it in \cite{mao1}. For more results on numerical methods, one can see e.g., \cite{MS1,LY,hut1,S,S1} and reference therein.

In this paper, we will present another explicit method for nonlinear SDEs. We call it the modified truncated EM (MTEM) method since it is motivated by and different from the truncated EM method introduced by Mao in \cite{mao}. Strong convergence rates of the MTEM method to the exact solutions to underlying SDEs are investigated under given condition. Results suggest that less conditions are needed to ensure the strong convergence for the MTEM method than the truncated EM method.

The organization of the paper is as the following. In Section 2, the MTEM method is developed, and main results are presented. In Section 3, some useful lemmas are presented to obtain the convergence theorems. In Section 4, convergence rates at fixed time $T$ are obtained. Then the convergence rates over the time interval $[0,T]$ will be proved with additional assumption on $g$ in Section 5. In Section 5. Then in Section 6, two examples will be presented to interpret the theory. We will conclude our paper in Section 7.

\section{The settings and main results}

Let $(\Omega,\mathscr{F},(\mathscr{F}_t)_{t\geq 0},P)$ be a complete
filtered probability space satisfying usual conditions. Consider the following stochastic differential
equations:

\begin{equation}\label{sde}dX(t)
=f(X(t))dt+g(X(t))dB_t,\ X_0=x_0\in\mathbb{R}^d
\end{equation}
where  $(B_t)_{t\geq0}$ is an $m$-dimensional standard
$\mathscr{F}_t$-Brownian motion,
$f:x\in\mathbb{R}^d\mapsto
f(x)\in\mathbb{R}^d$ and
$g:x\in\mathbb{R}^d\mapsto g(x)\in\mathbb{R}^d
\otimes\mathbb{R}^m$ are measurable functions.

Assume that the coefficients satisfy local Lipschitz condition, that is, for each
$R>0$ there is $L_R>0$ (depending on $R$) such that
\begin{equation}\label{local}|f(x)-f(\bar{x})|
\vee|g(x)-g(\bar{x})|\le L_R|x-\bar{x}|\end{equation} for all
$|x|\vee|\bar{x}|\le R$, where $|x|$ is the Euclidean norm for vector $x\in\mathbb{R}^d$ and $|A|=\sqrt{\textrm{trace}(A^TA)}$ is the trace norm for a matrix $A$.

It is obvious that $L_R$ is an increasing function with respect to $R$, we only need to consider that case that $L_R\uparrow\infty$ as $R\to\infty$ for simplicity. It is also well known that there is a unique strong solution (might explode at finite time) to equation (\ref{sde}) under local Lipschitz condition (\ref{local}) (Indeed, local Lipschitz condition could be relaxed to non-Lipschitz condition, see e.g.
\cite{LW}).

Choose $\Delta^*>0$ small enough and a strictly positive decreasing function $h:(0,\Delta^*]\to(0,\infty)$ such that
\begin{equation}\label{tiaoj}\lim_{\Delta\to0}h(\Delta)=\infty\ \textrm{and}\ \lim_{\Delta\to0}L_{h(\Delta)}^4\Delta=0.\end{equation}

\begin{Remark}\label{r1}
Such function $h$ always exists for any given Lipschitz coefficient $L_R$. Indeed, we can set $l(R)=\frac{1}{RL_R^4}$ and $h$ is the inverse function of $l$. Then $h$ is decreasing and $\lim_{\Delta\to0}h(\Delta)=\infty$ since $l$ is decreasing and $\lim_{R\to\infty}l(R)=0$. If we set $R=h(\Delta),$ then $L_{h(\Delta)}^4\Delta=L_R^4l(R)=\frac{1}{R}=\frac{1}{h(\Delta)}\to0$ as $\Delta\to0$. For example, let $L_R=2R.$ We define $h(x):=\frac{1}{\sqrt[5]{16x}}$. It is clear that $\lim_{\Delta\to0}h(\Delta)=\lim_{\Delta\to0}\frac{1}{\sqrt[5]{16\Delta}}=\infty$, and $L_{h(\Delta)}^4\Delta=\frac{1}{h(\Delta)}=\sqrt[5]{16\Delta}\to0.$ That is, (\ref{tiaoj}) holds for such defined $h.$
\end{Remark}

Motivated by Mao \cite{mao}, for any fixed $\Delta>0,$ we define the modified truncated function of $f$ as the following:
$$f_\Delta(x)=\left\{\begin{array}{ll} f(x),|x|\le h(\Delta),\\
\frac{|x|}{h(\Delta)} \cdot f(h(\Delta)\cdot\frac{x}{|x|}),|x|>h(\Delta).
\end{array}
\right.$$
$g_\Delta$ is defined in the same way as $f_\Delta$.

Notice that the modified truncated functions of $f$ and $g$ defined above are different from Mao \cite{mao} (where the truncated functions are bound for any fixed $\Delta$).

For the modified truncated function $f_\Delta$ and $g_\Delta$, we have the following

\begin{Lemma}
Suppose the local Lipschitz condition (\ref{local}) holds. Then for fixed $\Delta>0$ (small enough such that $f(0)\le h(\Delta)$ and $L_{h(\Delta)}\ge1$), the modified truncated functions $f_\Delta$ and $g_\Delta$ are global Lipschitz continuous with Lipschitz constant $4L_{h(\Delta)}.$ That is
\begin{equation}\label{global}|f_\Delta(x)-f_\Delta(\bar{x})|
\vee|g_\Delta(x)-g_\Delta(\bar{x})|\le 4L_{h(\Delta)}|x-\bar{x}|,\ \forall x,\bar{x}\in\mathbb{R}^d.\end{equation}
\end{Lemma}

\textbf{Proof}\ For any $x,\bar{x}\in\mathbb{R}^d$, there are three cases: $x,\bar{x}$ are both in the ball $B(h(\Delta))=\{x\in\mathbb{R}^d,|x|\le h(\Delta)\},$ $x,\bar{x}$ are both outside the ball $B(h(\Delta))$ and one is in the ball and the other is outside the ball.

\textbf{Case 1}. $x,\bar{x}\le h(\Delta)$. Then
$$|f_\Delta(x)-f_\Delta(\bar{x})|=|f(x)-f(\bar{x})|\le L_{h(\Delta)}|x-\bar{x}|\le 4L_{h(\Delta)}|x-\bar{x}|.$$

\textbf{Case 2}. $x,\bar{x}> h(\Delta)$. Since
$$|h(\Delta)\cdot\frac{x}{|x|}|=|h(\Delta)\cdot\frac{\bar{x}}{|\bar{x}|}|=h(\Delta),$$

then we have

$$\aligned |f_\Delta(x)-f_\Delta(\bar{x})|&=\left|\frac{|x|}{h(\Delta)} \cdot f\left(h(\Delta)\cdot\frac{x}{|x|}\right)-\frac{|\bar{x}|}{h(\Delta)} \cdot f\left(h(\Delta)\cdot\frac{\bar{x}}{|\bar{x}|}\right)\right|\\&
\le \frac{|x|}{h(\Delta)}\left|f\left(h(\Delta)\cdot\frac{x}{|x|}\right)-f\left(h(\Delta)\cdot\frac{\bar{x}}{|\bar{x}|}\right)\right|\\&
\quad+ \frac{\left||x|-|\bar{x}|\right|}{h(\Delta)}\left(\left|f\left(h(\Delta)\cdot\frac{\bar{x}}{|\bar{x}|}\right)-f(0)\right|+|f(0)|\right)\\&
\le \frac{|x|}{h(\Delta)}\cdot L_{h(\Delta)}\left|h(\Delta)\cdot\frac{x}{|x|}-h(\Delta)\cdot\frac{\bar{x}}{|\bar{x}|}\right|\\&
\quad+ \frac{\left||x|-|\bar{x}|\right|}{h(\Delta)}\left(L_{h(\Delta)}\left|h(\Delta)\cdot\frac{\bar{x}}{|\bar{x}|}\right|+|f(0)|\right)\\&
\le L_{h(\Delta)}\left|x-\frac{|x|}{|\bar{x}|}\bar{x}\right|+\frac{L_{h(\Delta)}h(\Delta)+|f(0)|}{h(\Delta)}|x-\bar{x}|\\&
\le L_{h(\Delta)}\left(\left|x-\bar{x}\right|+\left|\bar{x}-\frac{|x|}{|\bar{x}|}\bar{x}\right|\right)
+\frac{L_{h(\Delta)}h(\Delta)+|f(0)|}{h(\Delta)}|x-\bar{x}|\\&
\le (3L_{h(\Delta)}+1)|x-\bar{x}|\le 4L_{h(\Delta)}|x-\bar{x}|.\endaligned$$

\textbf{Case 3}. One is in the ball and the other is outside the ball. Without loss of generality, suppose that $|x|\le h(\Delta)<|\bar{x}|.$ Then we have
$$\aligned |f_\Delta(x)-f_\Delta(\bar{x})|&=\left|f(x)-\frac{|\bar{x}|}{h(\Delta)} \cdot f\left(h(\Delta)\cdot\frac{\bar{x}}{|\bar{x}|}\right)\right|\\&
\le \left|f(x)-f\left(h(\Delta)\cdot\frac{\bar{x}}{|\bar{x}|}\right)\right|
+\left|f\left(h(\Delta)\cdot\frac{\bar{x}}{|\bar{x}|}\right)\right|\left|1-\frac{|\bar{x}|}{h(\Delta)}\right|\\&
\le L_{h(\Delta)}\left|x-h(\Delta)\cdot\frac{\bar{x}}{|\bar{x}|}\right|+\frac{L_{h(\Delta)}h(\Delta)+|f(0)|}{h(\Delta)}|h(\Delta)-|\bar{x}||.
\endaligned$$

Since $|x|\le h(\Delta)<|\bar{x}|$, then $|h(\Delta)-|\bar{x}||=|\bar{x}|-h(\Delta)\le|\bar{x}|-|x|\le|x-\bar{x}|$, and $$\aligned\left|x-h(\Delta)\cdot\frac{\bar{x}}{|\bar{x}|}\right|&\le |x-\bar{x}|+\left|\bar{x}-h(\Delta)\cdot\frac{\bar{x}}{|\bar{x}|}\right|\\&
=|x-\bar{x}|+|h(\Delta)-|\bar{x}||\le2|x-\bar{x}|.\endaligned$$

Therefore,
$$\aligned|f_\Delta(x)-f_\Delta(\bar{x})|&\le \left(2L_{h(\Delta)}+\frac{L_{h(\Delta)}h(\Delta)+|f(0)|}{h(\Delta)}\right)|x-\bar{x}|
\\&\le (3L_{h(\Delta)}+1)|x-\bar{x}|\le 4L_{h(\Delta)}|x-\bar{x}|.\endaligned$$

We complete the proof. $\square$

Then we define the modified truncated EM (MTEM) method numerical solutions $X_k^\Delta\approx x(k\Delta)$ by setting $X_0^\Delta=x_0$ and
\begin{equation}\label{num}X_{k+1}^\Delta=X_k^\Delta+f_\Delta(X_k^\Delta)\Delta+g_\Delta(X_k^\Delta)\Delta B_k\end{equation}
for $k=0,1,2,\cdots,$ where $\Delta B_k=B((k+1)\Delta)-B(k\Delta)$ is the increment of the Brownian motion.

The two versions of the continuous-time MTEM solutions are defined as the following:

\begin{equation}\label{num1}\bar{x}_\Delta(t)=\sum_{k=0}^\infty X_k^\Delta1_{[k\Delta,(k+1)\Delta)}(t),\quad t\ge0,\end{equation}
and
\begin{equation}\label{num2}x_\Delta(t)=x_0+\int_0^tf_\Delta(\bar{x}_\Delta(s))ds+\int_0^tg_\Delta(\bar{x}_\Delta(s))dB(s),\quad t\ge0.\end{equation}

It is easy to see that $x_\Delta(k\Delta)=\bar{x}_\Delta(k\Delta)=X_k^\Delta$ for all $k\ge0.$

To study the strong convergence of MTEM (\ref{num}), let us consider the following conditions:

Suppose there exists $q\ge2$ and $H>0$ such that

\begin{equation}\label{lianxu}
\langle x-y,f(x)-f(y)\rangle+\frac{q-1}{2}|g(x)-g(y)|^2\le H|x-y|^2
\end{equation}
for all $x,y\in\mathbb{R}^d$ and there is a pair of constants $p>2$ and $K>0$ such that
\begin{equation}\label{zengzhang}
\langle x,f(x)\rangle+\frac{p-1}{2}|g(x)|^2\le K(1+|x|^2), x\in\mathbb{R}^d.
\end{equation}

Now we are ready to state our first result on the strong convergence rate for MTEM method at fixed time $T.$

\begin{Theorem}\label{conv}
Assume that (\ref{local}), (\ref{tiaoj}) and (\ref{lianxu}) hold, and (\ref{zengzhang}) holds for $2< p\le 6$. If there exists $2< q<p$ such that
\begin{equation}\label{tj}
h(\Delta)\ge(L^{2q}_{h(\Delta)}\Delta^{\frac{q}{2}})^{-\frac{1}{p-q}}
\end{equation}
holds for any $\Delta\le\Delta_0\ (\le\Delta^*)$, then the continuous-time MTEM methods satisfy
\begin{equation}\label{shou}\mathbb{E}|x(T)-x_\Delta(T)|^q\le C(T,q)L^{2q}_{h(\Delta)}\Delta^{q/2}\ \textrm{and}\ \mathbb{E}|x(T)-\bar{x}_\Delta(T)|^q\le C_{q,T}L^{2q}_{h(\Delta)}\Delta^{q/2}.\end{equation}
\end{Theorem}
\begin{Remark}
Notice that the set of $h(\Delta)$ such that (\ref{tj}) holds for small $\Delta$ is not empty. For example, let $L$ and $h$ be the same as in Remark \ref{r1}. Then we have $L^{4}_{h(\Delta)}\Delta=\frac{1}{h(\Delta)}.$ Thus (\ref{tj}) is equivalent to
$h(\Delta)\ge(h(\Delta))^{\frac{q}{2(p-q)}}$. Since $h(\Delta)\to\infty$ as $\Delta\to0,$ then if $\frac{q}{2(p-q)}\le1,$ (\ref{tj}) holds for $\Delta$ small enough.
\end{Remark}

For the convergence rates over the time interval $[0,T],$ we have to introduce an additional assumption.

Suppose there exist $r\ge2$ and $\bar{K}>0$ such that
\begin{equation}\label{gzeng}
|g(x)|^2\le\bar{K}(1+|x|^r), \forall x\in \mathbb{R}^d.
\end{equation}

Let us now present our second strong converge result for the continuous-time MTEM method. This time, the strong convergence rates over a time interval are obtained under given conditions.
\begin{Theorem}\label{conv1}
Assume that (\ref{local}), (\ref{tiaoj}), (\ref{lianxu}), (\ref{zengzhang}) and (\ref{gzeng}) hold. If there exist $2\le r<p\le6$ and $2< q\le p+2-r$ such that (\ref{tj}) holds for $\Delta$ small enough, then there exists $C$ (independent of $\Delta$) such that
\begin{equation}\label{shou1}\mathbb{E}\sup_{0\le t\le T}|x(t)-x_\Delta(t)|^q\le CL^{2q}_{h(\Delta)}\Delta^{q/2}\end{equation}
and if further $2<q\le4$, then
\begin{equation}\label{shou2}
\mathbb{E}\sup_{0\le t\le T}|x(t)-\bar{x}_\Delta(t)|^q\le CL^{q}_{h(\Delta)}\Delta^{q/2-1}.
\end{equation}
\end{Theorem}

\section{Some useful lemmas}

Firstly, we present a property of $f_\Delta$ and $g_\Delta$ similar to $f$ and $g$.

\begin{Lemma}\label{l1} For $\Delta$ small enough, condition (\ref{zengzhang}) implies
\begin{equation}\label{bijin}\langle x,f_\Delta(x)\rangle+\frac{p-1}{2}|g_\Delta(x)|^2\le 2K(1+|x|^2),\ \forall x\in\mathbb{R}^d.\end{equation}
\end{Lemma}

\textbf{Proof}\ First, suppose $|x|\le h(\Delta).$ Then we have $f_\Delta(x)=f(x)$ and $g_\Delta(x)=g(x)$. Thus
$$\langle x,f_\Delta(x)\rangle+\frac{p-1}{2}|g_\Delta(x)|^2=\langle x,f(x)\rangle+\frac{p-1}{2}|g(x)|^2\le K(1+|x|^2),\ \forall |x|\le h(\Delta).$$

If $|x|>h(\Delta)$, then $f_\Delta(x)=\frac{|x|}{h(\Delta)}f(h(\Delta)\frac{x}{|x|})$ and $g_\Delta(x)=\frac{|x|}{h(\Delta)}g(h(\Delta)\frac{x}{|x|})$. Notice that $\left|h(\Delta)\cdot \frac{x}{|x|}\right|=h(\Delta)$. Therefore
$$\aligned \langle x,f_\Delta(x)\rangle+\frac{p-1}{2}|g_\Delta(x)|^2&=\frac{|x|^2}{h^2(\Delta)}\left(\left\langle h(\Delta)\frac{x}{|x|},f\left(h(\Delta)\frac{x}{|x|}\right)\right\rangle
+\frac{p-1}{2}\left|g\left(h(\Delta)\frac{x}{|x|}\right)\right|^2\right)\\&
\le |x|^2\cdot K\left(\frac{1+h^2(\Delta)}{h^2(\Delta)}\right)\le 2K(1+|x|^2).\endaligned$$
We have used (\ref{zengzhang}) and (\ref{tiaoj}) in the last two inequalities, respectively.
Thus for any $x\in\mathbb{R}^d$, (\ref{bijin}) holds.
We complete the proof. $\square$

Now let us state the following two known results as lemmas (see \cite{mao,mao1}) for the the proof of Theorem \ref{conv}. First, we have

\begin{Lemma}\label{ju}
Under conditions (\ref{local}) and (\ref{zengzhang}), the SDE (\ref{sde}) has a unique global solution $x(t)$ and, moreover,
$$\sup_{0\le t\le T}\mathbb{E}|x(t)|^p<\infty,\ \forall T>0.$$
\end{Lemma}

\begin{Lemma}\label{ting}
Define the stopping time
$$\tau_R=\inf\{t\ge0,|x(t)|\ge R\},\ \inf\emptyset=\infty.$$
Suppose conditions (\ref{local}) and (\ref{zengzhang}) hold. Then
$$P(\tau_R\le T)\le\frac{C}{R^p}.$$
\end{Lemma}

As a similar result of Lemma \ref{ju}, we have the following moment property for the MTEM method (\ref{num2}).

\begin{Lemma}\label{temju}
Assume that conditions (\ref{local}), (\ref{tiaoj}) and (\ref{zengzhang}) hold for $0< p\le6$. Then there exist $0<\Delta_0\le\Delta^*$ and a constant $C(T,p)>0$ (independent of $\Delta$) such that for any $\Delta\in(0,\Delta_0]$, the modified TEM method (\ref{num}) satisfies
$$\sup_{0<\Delta\le\Delta_0}\sup_{0\le k\le[\frac{T}{\Delta}]}\mathbb{E}|X_k^\Delta|^p\le C(T,p)<\infty,\ \forall T>0.$$
\end{Lemma}

\textbf{Proof}\ By definition of (\ref{num}), for any $\Delta>0$ and any $0\le k\le [\frac{T}{\Delta}],$ we have

$$\aligned|X_{k+1}^\Delta|^2&=|X_{k}^\Delta|^2+2\langle X_{k}^\Delta,f_\Delta(X_k^\Delta)\Delta+g_\Delta(X_k^\Delta)\Delta B_k\rangle+|f_\Delta(X_k^\Delta)\Delta+g_\Delta(X_k^\Delta)\Delta B_k|^2\\&=|X_{k}^\Delta|^2+\xi_k,\endaligned$$
where
$$\xi_k:=2\langle X_{k}^\Delta,f_\Delta(X_k^\Delta)\Delta+g_\Delta(X_k^\Delta)\Delta B_k\rangle+|f_\Delta(X_k^\Delta)\Delta+g_\Delta(X_k^\Delta)\Delta B_k|^2.$$

Then we have
$$\aligned|X_{k+1}^\Delta|^p&=(|X_{k}^\Delta|^2+\xi_k)^\frac{p}{2}.\endaligned$$

Define $F(x)=(|X_{k}^\Delta|^2+x)^\frac{p}{2}$. Then for any $x\in\mathbb{R}^1$ and $p\in[4,6]$, by Taylor's expansion,
$$\aligned F(x)&=|X_{k}^\Delta|^p+\frac{p}{2}|X_{k}^\Delta|^{p-2}x+\frac{p(p-2)}{8}|X_{k}^\Delta|^{p-4}x^2\\&
\quad+\frac{p(p-2)(p-4)}{2^3\times3!}|X_{k}^\Delta|^{p-6}x^3
+\frac{p(p-2)(p-4)(p-6)}{2^4\times4!}|X_{k}^\Delta|^{p-8}\theta^4
\\&\le|X_{k}^\Delta|^p+\frac{p}{2}|X_{k}^\Delta|^{p-2}x+\frac{p(p-2)}{8}|X_{k}^\Delta|^{p-4}x^2\\&
\quad+\frac{p(p-2)(p-4)}{2^3\times3!}|X_{k}^\Delta|^{p-6}x^3,\endaligned$$
where $\theta$ in the first equation lies between 0 and $x$.

Then we have
$$\aligned\mathbb{E}(|X_{k+1}^\Delta|^p|\mathscr{F}_{k\Delta})&\le\mathbb{E}
\left(|X_{k}^\Delta|^p+\frac{p}{2}|X_{k}^\Delta|^{p-2}\xi_k+\frac{p(p-2)}{8}|X_{k}^\Delta|^{p-4}\xi_k^2\right.\\&
\qquad+\left.\frac{p(p-2)(p-4)}{2^3\times3!}|X_{k}^\Delta|^{p-6}\xi_k^3|\mathscr{F}_{k\Delta}\right).\endaligned$$

Now
$$\aligned\mathbb{E}(\xi_k|\mathscr{F}_{k\Delta})&=
\mathbb{E}(2\langle X_{k}^\Delta,f_\Delta(X_k^\Delta)\Delta+g_\Delta(X_k^\Delta)\Delta B_k\rangle+|f_\Delta(X_k^\Delta)\Delta+g_\Delta(X_k^\Delta)\Delta B_k|^2|\mathscr{F}_{k\Delta})\\&=2\langle X_{k}^\Delta, f_\Delta(X_k^\Delta)\rangle+|g_\Delta(X_k^\Delta)|^2)\Delta+|f_\Delta(X_k^\Delta)|^2\Delta^2.\endaligned$$

We have used the fact that $\mathbb{E}(\Delta B_k|\mathscr{F}_{k\Delta})=0$ and $\mathbb{E}(|\Delta B_k|^2|\mathscr{F}_{k\Delta})=\Delta$ in the above equation.

Since $f_\Delta$ satisfies the global Lipschitz condition (\ref{global}), then
$$\aligned\mathbb{E}(\xi_k|\mathscr{F}_{k\Delta})&\le (2\langle X_{k}^\Delta, f_\Delta(X_k^\Delta)\rangle+|g_\Delta(X_k^\Delta)|^2)\Delta
+(32L_{h(\Delta)}^2|X_{k}^\Delta|^2+2|f(0)|^2)\Delta^2.\endaligned$$

By (\ref{tiaoj}), we have
$$\frac{L^2_{h(\Delta)}\Delta^2}{\Delta}=L^2_{h(\Delta)}\Delta\to0\ \textrm{as}\ \Delta\to0.$$

That is $L^2_{h(\Delta)}\Delta^2=o(\Delta).$ Here and from now on, $o(\Delta)$ represents the higher order infinitesimal of $\Delta$ as $\Delta\to0.$ Therefore,

$$\aligned |X_{k}^\Delta|^{p-2}\mathbb{E}(\xi_k|\mathscr{F}_{k\Delta})&\le |X_{k}^\Delta|^{p-2}(2\langle X_{k}^\Delta, f_\Delta(X_k^\Delta)\rangle+|g_\Delta(X_k^\Delta)|^2)\Delta\\&\quad
+[64L^2_{h(\Delta)}|X_{k}^\Delta|^p+2|f(0)|^2|X_{k}^\Delta|^{p-2}]\Delta^2\\&
\le |X_{k}^\Delta|^{p-2}(2\langle X_{k}^\Delta, f_\Delta(X_k^\Delta)\rangle+|g_\Delta(X_k^\Delta)|^2)\Delta+o(\Delta)|X_{k}^\Delta|^p+o(\Delta).\endaligned$$
We have used the fact that $|x|^{i}\le 1+|x|^j, \forall 0<i<j, x\in\mathbb{R}^1.$

Similarly, by (\ref{tiaoj}), we have
$$\aligned\mathbb{E}(\xi^2_k|\mathscr{F}_{k\Delta})&=
\mathbb{E}((2\langle X_{k}^\Delta, g_\Delta(X_k^\Delta)\Delta B_k\rangle+B)^2|\mathscr{F}_{k\Delta})\\&\le
\mathbb{E}(4|X_{k}^\Delta|^2 |g_\Delta(X_k^\Delta)|^2|\Delta B_k|^2+B^2+4B\langle X_{k}^\Delta, g_\Delta(X_k^\Delta)\Delta B_k\rangle|\mathscr{F}_{k\Delta}),\endaligned$$
where $$B:=2\langle X_k^\Delta, f_\Delta(X_k^\Delta)\rangle\Delta+|f_\Delta(X_k^\Delta)|^2\Delta^2+|g_\Delta(X_k^\Delta)\Delta B_k|^2+2\langle f_\Delta(X_k^\Delta),g_\Delta(X_k^\Delta)\Delta B_k\rangle\Delta.$$

According to (\ref{tiaoj}) and the fact that $|x|^{i}\le 1+|x|^j, \forall 0<i<j, x\in\mathbb{R}^1$ again, it follows that
$$\aligned\mathbb{E}(B\langle X_{k}^\Delta, g_\Delta(X_k^\Delta)\Delta B_k\rangle|\mathscr{F}_{k\Delta})&\le2|X_{k}^\Delta| |f_\Delta(X_k^\Delta)||g_\Delta(X_k^\Delta)|^2\Delta^2\\&
\le 2|X_k^\Delta|(4L_{h(\Delta)}|X_{k}^\Delta|+|f(0)|)\\&\quad\times 2((4L_{h(\Delta)})^2|X_{k}^\Delta|^2+|g(0)|^2)\Delta^2\\&
\le 4[(4L_{h(\Delta)})^3|X_{k}^\Delta|^4+(4L_{h(\Delta)})^2|f(0)||X_{k}^\Delta|^3\\&
\quad+(4L_{h(\Delta)})|g(0)|^2|X_{k}^\Delta|^2+|f(0)||g(0)|^2|X_{k}^\Delta|]\Delta^2\\&
\le |X_{k}^\Delta|^4\cdot o(\Delta)+ o(\Delta),\endaligned$$
and in the same way
$$\aligned\mathbb{E}(B^2|\mathscr{F}_{k\Delta})&\le 2(4|X_k^\Delta|^2 |f_\Delta(X_k^\Delta)|^2\Delta^2+|f_\Delta(X_k^\Delta)|^4\Delta^4\\&\qquad
+3|g_\Delta(X_k^\Delta)|^4\Delta^2+4|f_\Delta(X_k^\Delta)|^2 |g_\Delta(X_k^\Delta)|^2\Delta^3)\\&
\le |X_{k}^\Delta|^4\cdot o(\Delta)+ o(\Delta).\endaligned$$

Then
$$\aligned |X_{k}^\Delta|^{p-4}\mathbb{E}(\xi^2_k|\mathscr{F}_{k\Delta})&\le |X_{k}^\Delta|^{p-2}\cdot4 |g_\Delta(X_k^\Delta)|^2\Delta+|X_{k}^\Delta|^p\cdot o(\Delta)+o(\Delta).\endaligned$$

Moreover, we can use the same method to derive that
$$\aligned |X_{k}^\Delta|^{p-6}\mathbb{E}(\xi^3_k|\mathscr{F}_{k\Delta})&\le |X_{k}^\Delta|^p\cdot o(\Delta)+o(\Delta).\endaligned$$

Therefore,

$$\aligned\mathbb{E}(|X_{k+1}^\Delta|^p|\mathscr{F}_{k\Delta})&\le|X_{k}^\Delta|^p+\frac{p}{2}|X_{k}^\Delta|^{p-2}\Big(2\langle X_{k}^\Delta, f_\Delta(X_k^\Delta)\rangle+(p-1)|g_\Delta(X_k^\Delta)|^2)\Big)\Delta
\\&\quad+|X_{k}^\Delta|^p\cdot o(\Delta)+o(\Delta)\\&
=|X_{k}^\Delta|^p+p|X_{k}^\Delta|^{p-2}\left(\langle X_{k}^\Delta, f_\Delta(X_k^\Delta)\rangle+\frac{p-1}{2}|g_\Delta(X_k^\Delta)|^2\right)\Delta\\&
\quad+|X_{k}^\Delta|^p\cdot o(\Delta)+o(\Delta).\endaligned$$

Then for any $0<\varepsilon(<1)$, we can choose $\Delta_0$ small enough such that for any $\Delta\le\Delta_0$, $o(\Delta)\le\varepsilon\Delta$. Now by condition (\ref{zengzhang}) and Lemma (\ref{l1}), we have
$$\aligned\mathbb{E}(|X_{k+1}^\Delta|^p|\mathscr{F}_{k\Delta})&\le |X_{k}^\Delta|^p+p|X_{k}^\Delta|^{p-2}\cdot 2K(1+|X_{k}^\Delta|^2)\Delta+|X_{k}^\Delta|^p\cdot \varepsilon\Delta+\varepsilon\Delta\\&
\le|X_{k}^\Delta|^p+2pK(|X_{k}^\Delta|^{p}+1)\Delta+2pK|X_{k}^\Delta|^{p}\Delta+|X_{k}^\Delta|^p\cdot \varepsilon\Delta+\varepsilon\Delta\\&
\le |X_{k}^\Delta|^p(1+(4pK+\varepsilon)\Delta)+(2pK+\varepsilon)\Delta.
\endaligned$$

Taking expectation on both sides, it follows that
$$\mathbb{E}(|X_{k+1}^\Delta|^p)\le(1+(4pK+\varepsilon)\Delta)\mathbb{E}(|X_{k}^\Delta|^p)+(2pK+\varepsilon)\Delta.$$

By induction, we have
$$\aligned\mathbb{E}(|X_{k}^\Delta|^p)&\le(1+(4pK+\varepsilon)\Delta)^k|x_0|^p+(2pK+\varepsilon)\Delta\sum_{i=0}^{k-1}(1+(4pK+\varepsilon)\Delta)^i
\\&\le e^{(4pK+\varepsilon)k\Delta}|x_0|^p+(2pK+\varepsilon)\Delta\frac{e^{(4pK+\varepsilon)k\Delta}-1}{(4pK+\varepsilon)\Delta}\\&
\le e^{(4pK+\varepsilon)T}(|x_0|^p+1).\endaligned$$

Set $C(T,p)=e^{(4pK+\varepsilon)T}(|x_0|^p+1)$. We have proved the conclusion for $4\le p\le6.$ For $0<p<4,$ by H\"{o}der's inequality, it follows that $$\sup_{0<\Delta\le\Delta_0}\sup_{0\le k\le[\frac{T}{\Delta}]}\mathbb{E}|X_k^\Delta|^p\le [\sup_{0<\Delta\le\Delta_0}\sup_{0\le k\le[\frac{T}{\Delta}]}\mathbb{E}(|X_k^\Delta|^4)]^\frac{p}{4}<C(p,T)<\infty.$$
This completes the proof. $\square$

Now let us present a lemma which shows that $x_\Delta(t)$ and $\bar{x}_\Delta(t)$ are close to each other in the sense of $L^p.$

\begin{Lemma}\label{close} Assume that (\ref{local}), (\ref{tiaoj}) and (\ref{zengzhang}) hold for $0< p\le6.$ For any $\Delta\in(0,\Delta^*)$, there exists $C(p,T)>0$ (independent of $\Delta$) such that
\begin{equation}\sup_{0\le t\le T}\mathbb{E}|x_\Delta(t)-\bar{x}_\Delta(t)|^p\le C(p,T)L_{h(\Delta)}^p\Delta^\frac{p}{2}.\end{equation}
\end{Lemma}

\textbf{Proof}\ For any fixed $t\le T,$ there exists $0\le k\le[\frac{T}{\Delta}]$ such that $k\Delta\le t<(k+1)\Delta.$ Thus
$$x_\Delta(t)-\bar{x}_\Delta(t)=x_\Delta(t)-X_k^\Delta=f_\Delta(X_k^\Delta)(t-k\Delta)+g_\Delta(X_k^\Delta)(B(t)-B(k\Delta)).$$

So we have
$$\mathbb{E}|x_\Delta(t)-\bar{x}_\Delta(t)|^p\le C_p(\Delta^p\mathbb{E}|f_\Delta(X_k^\Delta)|^p+\mathbb{E}(|g_\Delta(X_k^\Delta)|^p\mathbb{E}(|B(t)-B(k\Delta)|^p|\mathscr{F}_{k\Delta})).$$

Since $f_\Delta$ and $g_\Delta$ satisfy the global Lipschitz condition (\ref{global}), and notice that $B(t)-B(k\Delta)$ is independent of $\mathscr{F}_{k\Delta},$ then
$$\aligned\mathbb{E}|x_\Delta(t)-\bar{x}_\Delta(t)|^p&\le C_p\big(\Delta^p\mathbb{E}(4L_{h(\Delta)}|X_k^\Delta|+|f(0)|)^p\\&\quad
+\mathbb{E}(4L_{h(\Delta)}|X_k^\Delta|+|g(0)|)^p\Delta^\frac{p}{2}\big)\\&
\le C_p[4^pL_{h(\Delta)}^p\Delta^p(\mathbb{E}(|X_k^\Delta|^p)+|f(0)|^p)\\&\quad
+4^pL_{h(\Delta)}^p\Delta^\frac{p}{2}(\mathbb{E}(|X_k^\Delta|^p)+|g(0)|^p)].\endaligned$$

Therefore,
$$\sup_{0\le t\le T}\mathbb{E}|x_\Delta(t)-\bar{x}_\Delta(t)|^p\le C_pL^p_{h(\Delta)}\Delta^\frac{p}{2}(\sup_{0\le k\le[\frac{T}{\Delta}]}\mathbb{E}(|X_k^\Delta|^p)+|g(0)|^p),$$
where $C_p$ is a positive constant (independent of $\Delta$) which might change the value from line to line.
Then by Lemma \ref{temju}, we have
$$\sup_{0\le t\le T}\mathbb{E}|x_\Delta(t)-\bar{x}_\Delta(t)|^p\le C(p,T)L^p_{h(\Delta)}\Delta^\frac{p}{2}.$$

We then complete the proof. $\square$

For continuous-time MTEM method (\ref{num2}), we also have
\begin{Lemma}\label{temju1}
Assume that conditions (\ref{local}), (\ref{tiaoj}) and (\ref{zengzhang}) hold for $0< p\le6$. Then for any $\Delta\in(0,\Delta^*)$, there exists a constant $C(p,T)>0$ and $\Delta^*$ small enough such that the modified TEM (\ref{num}) satisfies
\begin{equation}\label{quejie}\sup_{0<\Delta\le\Delta^*}\sup_{0\le t\le T}\mathbb{E}|x_\Delta(t)|^p<C(p,T)<\infty,\ \forall T>0.\end{equation}
\end{Lemma}
\textbf{Proof}\ Notice that for any $0\le t\le T,$ there exists $k\le[\frac{T}{\Delta}]$ such that $k\Delta\le t<(k+1)\Delta.$ Thus $$\mathbb{E}|x_\Delta(t)|^p\le C_p(\mathbb{E}|x_\Delta(t)-\bar{x}_\Delta(t)|^p+\mathbb{E}|X_k^\Delta|^p).$$
Then (\ref{quejie}) follows directly by Lemma \ref{temju} and \ref{close}. $\square$

As a similar result of Lemma \ref{ting}, we have
\begin{Lemma}\label{ting1}
Define the stopping time
$$\rho_{\Delta,R}=\inf\{t\ge0,|x_\Delta(t)|\ge R\}.$$
Suppose conditions (\ref{local}), (\ref{tiaoj}) and (\ref{zengzhang}) hold for $0<p\le6$. Then for any $R>|x_0|$ and $\Delta\in(0,\Delta^*)$ ($\Delta^*$ small enough), we have
$$P(\rho_{\Delta,R}\le T)\le\frac{C}{R^p}.$$
\end{Lemma}

\textbf{Proof}\ We simply write $\rho_{\Delta,R}=\rho.$ By It\^o formula and Lemma \ref{l1}, for any $0\le t\le T,$
$$\aligned\mathbb{E}(|x_\Delta(t\wedge\rho)|^p)&\le|x_0|^p+\frac{p}{2}\mathbb{E}\int_0^{t\wedge\rho}|x_\Delta(s)|^{p-2}
(2x_\Delta(s)f_\Delta(\bar{x}_\Delta(s))+(p-1)|g_\Delta(\bar{x}_\Delta(s))|^2)ds\\&
=|x_0|^p+\frac{p}{2}\mathbb{E}\int_0^{t\wedge\rho}|x_\Delta(s)|^{p-2}
(2\bar{x}_\Delta(s)f_\Delta(\bar{x}_\Delta(s))+(p-1)|g_\Delta(\bar{x}_\Delta(s))|^2)ds\\&
\quad+p\mathbb{E}\int_0^{t\wedge\rho}|x_\Delta(s)|^{p-2}
(x_\Delta(s)-\bar{x}_\Delta(s))f_\Delta(\bar{x}_\Delta(s))ds\\&
\le|x_0|^p+\frac{p}{2}\mathbb{E}\int_0^{t\wedge\rho}|x_\Delta(s)|^{p-2}\cdot 2K(1+|\bar{x}_\Delta(s)|^2)ds\\&
\quad+p\mathbb{E}\int_0^{t\wedge\rho}|x_\Delta(s)|^{p-2}
|x_\Delta(s)-\bar{x}_\Delta(s)|\cdot (4L_{h(\Delta)}|\bar{x}_\Delta(s)|+|f(0)|)ds\\&
\le|x_0|^p+3pK\mathbb{E}\int_0^{t\wedge\rho}(1+|x_\Delta(s)|^{p})ds\\&\quad+2pK\mathbb{E}\int_0^{t\wedge\rho}|x_\Delta(s)|^{p-2}
|\bar{x}_\Delta(s)-x_\Delta(s)|^2ds\\&
\quad+p\mathbb{E}\int_0^{t\wedge\rho}4L_{h(\Delta)}|x_\Delta(s)-\bar{x}_\Delta(s)||x_\Delta(s)|^{p-2}
|\bar{x}_\Delta(s)|ds
\\&\quad+p|f(0)|\mathbb{E}\int_0^{t\wedge\rho}|x_\Delta(s)|^{p-2}|x_\Delta(s)-\bar{x}_\Delta(s)|ds\\&
=:|x_0|^p+J_1+J_2+J_3+J_4.\endaligned$$
Where $J_i$ denotes the $i$-th expectation of integral in the above expressions, $i=1,2,3,4$.

Notice that by Young's inequality, we have
$$\aligned J_2&\le 3pK\mathbb{E}\int_0^{t\wedge\rho}(\frac{p-2}{p}|x_\Delta(s)|^{p}+\frac{2}{p}|\bar{x}_\Delta(s)-x_\Delta(s)|^p)ds\\&
=3K(p-2)\mathbb{E}\int_0^{t\wedge\rho}(1+|x_\Delta(s)|^{p})ds+6K\mathbb{E}\int_0^{t\wedge\rho}|\bar{x}_\Delta(s)-x_\Delta(s)|^pds.\endaligned$$

Moreover,
$$\aligned J_4&\le p|f(0)|\mathbb{E}\int_0^{t\wedge\rho}|x_\Delta(s)|^{p-2}(1+|x_\Delta(s)-\bar{x}_\Delta(s)|^2)ds\\&
\le p|f(0)|\mathbb{E}\int_0^{t\wedge\rho}(1+|x_\Delta(s)|^p)ds+\frac{|f(0)|}{2K}J_2\\&
\le(\frac{5}{2}p-1)|f(0)|\mathbb{E}\int_0^{t\wedge\rho}(1+|x_\Delta(s)|^p)ds
+3|f(0)|\mathbb{E}\int_0^{t\wedge\rho}|\bar{x}_\Delta(s)-x_\Delta(s)|^pds.\endaligned$$

On the other hand, by using Young's inequality two times, we have
$$\aligned &\quad 4L_{h(\Delta)}|x_\Delta(s)-\bar{x}_\Delta(s)||x_\Delta(s)|^{p-2}
|\bar{x}_\Delta(s)|\\&\le\frac{1}{p}4^pL_{h(\Delta)}^p|x_\Delta(s)-\bar{x}_\Delta(s)|^p+\frac{1}{q}|\bar{x}_\Delta(s)|^q|x_\Delta(s)|^{q(p-2)}
\\&\le\frac{1}{p}4^pL_{h(\Delta)}^p|x_\Delta(s)-\bar{x}_\Delta(s)|^p
+\frac{1}{q}\left(\frac{q}{p}|\bar{x}_\Delta(s)|^p+(1-\frac{q}{p})|x_\Delta(s)|^{q(p-2)\frac{p}{p-q}}\right),\endaligned$$
where $q>1$ is a constant such that
$$\frac{1}{p}+\frac{1}{q}=1.$$

Thus $q(p-2)\frac{p}{p-q}=p.$ So we have
$$\aligned J_3&\le 4^pL_{h(\Delta)}^p\mathbb{E}\int_0^{t\wedge\rho}|x_\Delta(s)-\bar{x}_\Delta(s)|^pds
+\mathbb{E}\int_0^{t\wedge\rho}|\bar{x}_\Delta(s)|^pds+\frac{p-q}{q}\mathbb{E}\int_0^{t\wedge\rho}|x_\Delta(s)|^pds\\&
\le 4^pL_{h(\Delta)}^p\mathbb{E}\int_0^{T}|x_\Delta(s)-\bar{x}_\Delta(s)|^pds
+\mathbb{E}\int_0^{T}|\bar{x}_\Delta(s)|^pds+\frac{p-q}{q}\mathbb{E}\int_0^{t}|x_\Delta(s\wedge\rho)|^pds.\endaligned$$

Since $\sup_{0\le t\le T}\mathbb{E}|\bar{x}_\Delta(t)|^p=\sup_{0\le k\le [\frac{T}{\Delta}]}\mathbb{E}|X_k^\Delta|^p,$ then by Lemma \ref{temju} and \ref{close}, we have
$$J_3\le
C_{T,p}L^{2p}_{h(\Delta)}\Delta^{\frac{p}{2}}+C_{T,p}+\frac{p-q}{q}\int_0^{t}(1+\mathbb{E}|x_\Delta(s\wedge\rho)|^p)ds.$$

Therefore,
$$\mathbb{E}(|x_\Delta(t\wedge\rho)|^p)\le C_1+C_2\int_0^{t}(1+\mathbb{E}|x_\Delta(s\wedge\rho)|^p)ds.$$

Gronwall's lemma yields that
$$\mathbb{E}(|x_\Delta(T\wedge\rho)|^p)\le (C_1+1)e^{C_2T}-1<\infty.$$

This implies the required assertion easily. $\square$

\section{Convergence rate at fixed time $T$}

Let us first present a lemma which will play a key role in the proof of the convergence rate.

\begin{Lemma}\label{jubu}
Suppose (\ref{local}), (\ref{tiaoj}), (\ref{lianxu}) and (\ref{zengzhang}) hold for $2<q\le p\le6.$
Set $$\theta_{\Delta,R}=\tau_R\wedge\rho_{\Delta,R}\quad and\quad e_\Delta(t)=x(t)-x_\Delta(t)\ for\ t\ge0.$$
Then for any $\Delta\in(0,\Delta^*)$ and any $R\le h(\Delta)$, there exists $C_{q,T}>0$ (independent of $\Delta$) such that
$$\sup_{0\le t\le T}\mathbb{E}(|e_\Delta(t\wedge\theta_{\Delta,R})|^q)\le C_{q,T}L^{2q}_{h(\Delta)}\Delta^{\frac{q}{2}}.$$
\end{Lemma}

\textbf{Proof}\ Denote $\theta=\theta_{\Delta,R}$ for simplicity. By It\^o formula, we have

$$\aligned\mathbb{E}(|e_\Delta(t\wedge\theta)|^q)&\le\frac{q}{2}\mathbb{E}\int_0^{t\wedge\theta}|e_\Delta(s)|^{q-2}
[2\langle e_\Delta(s),f(x(s))-f_\Delta(\bar{x}_\Delta(s))\rangle\\&\qquad+(q-1)|g(x(s))-g_\Delta(\bar{x}_\Delta(s))|^2]ds.\endaligned$$

Since $\forall 0\le s\le t\wedge\theta,$ $|\bar{x}_\Delta(s)|\le R\le h(\Delta),$ then
$$f_\Delta(\bar{x}_\Delta(s))=f(\bar{x}_\Delta(s)),\quad g_\Delta(\bar{x}_\Delta(s))=g(\bar{x}_\Delta(s)),\quad \forall 0\le s\le t\wedge\theta.$$

Therefore,
$$\aligned\mathbb{E}(|e_\Delta(t\wedge\theta)|^q)&\le\frac{q}{2}\mathbb{E}\int_0^{t\wedge\theta}|e_\Delta(s)|^{q-2}
[2\langle e_\Delta(s),f(x(s))-f(\bar{x}_\Delta(s))\rangle\\&\qquad+(q-1)|g(x(s))-g(\bar{x}_\Delta(s))|^2]ds\\&
=q\mathbb{E}\int_0^{t\wedge\theta}|e_\Delta(s)|^{q-2}
[\langle x(s)-\bar{x}_\Delta(s),f(x(s))-f(\bar{x}_\Delta(s))\rangle\\&\qquad+\frac{q-1}{2}|g(x(s))-g(\bar{x}_\Delta(s))|^2]ds\\&
\quad+q\mathbb{E}\int_0^{t\wedge\theta}|e_\Delta(s)|^{q-2}
\langle\bar{x}_\Delta(s)-x_\Delta(s),f(x(s))-f(\bar{x}_\Delta(s))\rangle ds\\&
\le qH\mathbb{E}\int_0^{t\wedge\theta}|e_\Delta(s)|^{q-2}|x(s)-\bar{x}_\Delta(s)|^2ds\\&\quad
+q\mathbb{E}\int_0^{t\wedge\theta}|e_\Delta(s)|^{q-2}
|\bar{x}_\Delta(s)-x_\Delta(s)|\cdot4L_{h(\Delta)}|(x(s)-\bar{x}_\Delta(s)|ds\\&
\le 2qH\mathbb{E}\int_0^{t\wedge\theta}(|e_\Delta(s)|^{q}+|e_\Delta(s)|^{q-2}|x_\Delta(s)-\bar{x}_\Delta(s)|^2)ds\\&\quad
+q\mathbb{E}\int_0^{t\wedge\theta}|e_\Delta(s)|^{q-2}
|\bar{x}_\Delta(s)-x_\Delta(s)|\cdot4L_{h(\Delta)}|x(s)-\bar{x}_\Delta(s)|ds.\endaligned$$

By Young's inequality, we have
$$\aligned\mathbb{E}(|e_\Delta(t\wedge\theta)|^q)&\le (4q-4)H\int_0^{t}\mathbb{E}(|e_\Delta(s\wedge\theta)|^{q})ds+4H\int_0^{t}\mathbb{E}(|\bar{x}_\Delta(s)-x_\Delta(s)|^q)ds\\&
\quad+(q-2)\int_0^{t}\mathbb{E}(|e_\Delta(s\wedge\theta)|^{q})ds\\&\quad+2\mathbb{E}\int_0^{t\wedge\theta}
|\bar{x}_\Delta(s)-x_\Delta(s)|^\frac{q}{2}(4L_{h(\Delta)})^\frac{q}{2}|x(s)-\bar{x}_\Delta(s)|^\frac{q}{2}ds.\endaligned$$

Since
$$|x(s)-\bar{x}_\Delta(s)|^\frac{q}{2}\le C_{\frac{q}{2}}(|x_\Delta(s)-\bar{x}_\Delta(s)|^\frac{q}{2}+|x(s)-x_\Delta(s)|^\frac{q}{2}),$$
then
$$\aligned&\quad\mathbb{E}\int_0^{t\wedge\theta}
|\bar{x}_\Delta(s)-x_\Delta(s)|^\frac{q}{2}(4L_{h(\Delta)})^\frac{q}{2}|x(s)-\bar{x}_\Delta(s)|^\frac{q}{2}ds\\&
\le C_{\frac{q}{2}}[(L^\frac{q}{2}_{h(\Delta)}+L^q_{h(\Delta)})\int_0^{T}
\mathbb{E}|\bar{x}_\Delta(s)-x_\Delta(s)|^qds+\int_0^{T}
\mathbb{E}|e_\Delta(s\wedge\theta)|^qds].\endaligned$$

Notice that $q<p,$ then by Lemma \ref{close} and H\"{o}lder inequality,
$$\aligned\mathbb{E}(|e_\Delta(t\wedge\theta)|^q)&\le (4H(q-1)+q-2+C_\frac{q}{2})\int_0^{t}\mathbb{E}(|e_\Delta(s\wedge\theta)|^{q})ds+4HTL^q_{h(\Delta)}\Delta^\frac{q}{2}\\&
\quad+TC_{\frac{q}{2}}(L^\frac{3q}{2}_{h(\Delta)}+L^{2q}_{h(\Delta)})\Delta^\frac{q}{2}.\endaligned$$

By Gronwall's lemma, we have
$$\mathbb{E}(|e_\Delta(t\wedge\theta)|^q)\le C_{q,T}L^{2q}_{h(\Delta)}\Delta^\frac{q}{2},\quad \forall 0\le t\le T.$$

This completes the proof. $\square$

Now we are ready to prove Theorem \ref{conv}.

\textbf{Proof of Theorem \ref{conv}}\ Let $\tau_R, \rho_{\Delta,R}, \theta_{\Delta,R}$ and $e_\Delta(t)$ be the same as before.
Then by Young's inequality, we have that for any $\delta>0,$
$$\aligned\mathbb{E}(|e_\Delta(T)|^q)&\le \mathbb{E}(|e_\Delta(T)|^q1_{\{\theta_{\Delta,R}>T\}})+ \mathbb{E}(|e_\Delta(T)|^q1_{\{\theta_{\Delta,R}\le T\}})\\&\le\mathbb{E}(|e_\Delta(T)|^q1_{\{\theta_{\Delta,R}>T\}})+ \frac{q\delta}{p}\mathbb{E}(|e_\Delta(T)|^p)+\frac{p-q}{p\delta^{q/(p-q)}}P(\theta_{\Delta,R}\le T)\\&
\le\mathbb{E}(|e_\Delta(T\wedge\theta)|^q)+ \frac{q\delta C}{p}\left(\mathbb{E}(|x_\Delta(T)|^p)+\mathbb{E}(|x(T)|^p)\right)+\frac{p-q}{p\delta^{q/(p-q)}}P(\theta_{\Delta,R}\le T).\endaligned$$

By Lemma \ref{ju} and \ref{temju1}, we have
$$\mathbb{E}(|x_\Delta(T)|^p)+\mathbb{E}(|x(T)|^p\le C,$$
while by Lemma \ref{ting} and \ref{ting1},
$$P(\theta_{\Delta,R}\le T)\le P(\tau_R\le T)+P(\rho_{\Delta,R}\le T)\le\frac{C}{R^p}.$$

Thus,
$$\mathbb{E}(|e_\Delta(T)|^q)\le\mathbb{E}(|e_\Delta(T\wedge\theta)|^q)+\frac{qC\delta}{p}+\frac{C(p-q)}{pR^p\delta^{q/(p-q)}}$$
holds for any $\Delta\in(0,\Delta^*), R>|x_0|$ and $\delta>0$. Then we can choose $\delta=L^{2q}_{h(\Delta)}\Delta^\frac{q}{2}$ and $R=(L^{2q}_{h(\Delta)}\Delta^\frac{q}{2})^{-\frac{1}{p-q}}$ to get
$$\mathbb{E}(|e_\Delta(T)|^q)\le\mathbb{E}(|e_\Delta(T\wedge\theta)|^q)+CL^{2q}_{h(\Delta)}\Delta^\frac{q}{2}.$$

But by condition (\ref{tj}), we have
$$h(\Delta)\ge (L^{2q}_{h(\Delta)}\Delta^\frac{q}{2})^{-\frac{1}{p-q}}=R.$$

Then by Lemma \ref{jubu},
$$\mathbb{E}(|e_\Delta(T)|^q)\le CL^{2q}_{h(\Delta)}\Delta^\frac{q}{2}.$$

The second inequality there follows easily from the first one and Lemma \ref{close}. $\square$

\section{Convergence rates over the time interval $[0,T]$}

First of all, let us cite a Lemma from \cite{mao}.

\begin{Lemma}\label{zuida}
Let (\ref{local}), (\ref{zengzhang}) and (\ref{gzeng}) hold and assume that $p>r.$ Set $\bar{p}=2+p-r.$ Then
$$\mathbb{E}(\sup_{0\le t\le T}|x(t)|^{\bar{p}})<C,\ \forall T>0.$$
\end{Lemma}

\begin{Remark}
When $r=2,$ this result follows from Theorem 1.7 in \cite{LW} since our condition (\ref{zengzhang}) and (\ref{gzeng}) are stronger than that of (1.7) in Theorem 1.7 in \cite{LW} (notice that $p>2$). Indeed, they proved that $\mathbb{E}(\sup_{0\le t\le T}|x(t)|^p)<C.$
\end{Remark}

For the discontinuous and continuous-time MTEM methods (\ref{num}) and (\ref{num2}), we have
\begin{Lemma}\label{zuida1}
Let (\ref{local}), (\ref{tiaoj}), (\ref{zengzhang}) and (\ref{gzeng}) hold and assume that $6\ge p>r\ge2.$ Set $\bar{p}=2+p-r.$ Then
\begin{equation}\label{zuida2}\sup_{0<\Delta\le\Delta^*}\mathbb{E}(\sup_{0\le t\le T}|x_\Delta(t)|^{\bar{p}})<C,\ \forall T>0,\end{equation}
and therefore,
\begin{equation}\label{zuida3}\sup_{0<\Delta\le\Delta^*}\mathbb{E}(\sup_{0\le k\le [\frac{T}{\Delta}]}|X_k^\Delta|^{\bar{p}})<C,\ \forall T>0,\end{equation}
\end{Lemma}

\textbf{Proof}\ For any $\Delta\in(0,\Delta^*],$ by It\^o formula and Lemma \ref{l1}, we have
$$\aligned\sup_{0\le t\le T}|x_\Delta(t)|^{\bar{p}}&\le|x_0|^{\bar{p}}+2K\bar{p}\int_0^T|x_\Delta(s)|^{\bar{p}-2}
(1+|\bar{x}_\Delta(s)|^2)ds\\&
\quad+\bar{p}\int_0^T|x_\Delta(s)|^{\bar{p}-2}|x_\Delta(s)-\bar{x}_\Delta(s)||f_\Delta(\bar{x}_\Delta(s))|ds\\&
\quad+\bar{p}\sup_{0\le t\le T}\left|\int_0^t|x_\Delta(s)|^{\bar{p}-2}\langle x_\Delta(s),g_\Delta(\bar{x}_\Delta(s))dB(s)\rangle\right|.\endaligned$$

Notice that by Lemma \ref{temju1} and \ref{temju},
$$\mathbb{E}\int_0^T|x_\Delta(s)|^{\bar{p}-2}
(1+|\bar{x}_\Delta(s)|^2)ds\le C.$$
Moreover, since for $\Delta>0$ small enough,
$$|f_\Delta(x)|\le 4L_{h(\Delta)}|x|+|f(0)|,$$
then as in the proof of Lemma \ref{ting1}, we have
$$\mathbb{E}\int_0^T|x_\Delta(s)|^{\bar{p}-2}|x_\Delta(s)-\bar{x}_\Delta(s)||f_\Delta(\bar{x}_\Delta(s))|ds\le C.$$

So
$$\aligned\sup_{0\le t\le T}|x_\Delta(t)|^{\bar{p}}&\le C+\bar{p}\sup_{0\le t\le T}\left|\int_0^t|x_\Delta(s)|^{\bar{p}-2}\langle x_\Delta(s),g_\Delta(\bar{x}_\Delta(s))dB(s)\rangle\right|,\endaligned$$
where $C$ is a constant (independent of $\Delta$).
By the Burkholder-Davis-Gundy inequality (see e.g. \cite{Ikeda}) and (\ref{gzeng}), we have
$$\aligned\mathbb{E}\sup_{0\le t\le T}|x_\Delta(t)|^{\bar{p}}&\le C+4\sqrt{2}\bar{p}\mathbb{E}\left|\int_0^T|x_\Delta(s)|^{2\bar{p}-2}\bar{K}(1+|\bar{x}_\Delta(s)|^r)ds\right|^\frac{1}{2}\\&
\le C+\mathbb{E}\left|\sup_{0\le t\le T}|x_\Delta(t)|^{\bar{p}}\cdot\left(32\bar{p}^2\bar{K}\int_0^T|x_\Delta(s)|^{\bar{p}-2}
(1+|\bar{x}_\Delta(s)|^r)ds\right)\right|^\frac{1}{2}\\&
\le C+\frac{1}{2}\mathbb{E}\sup_{0\le t\le T}|x_\Delta(t)|^{\bar{p}}+16\bar{K}\bar{p}^2\mathbb{E}\int_0^T|x_\Delta(s)|^{\bar{p}-2}(1+|\bar{x}_\Delta(s)|^r)ds.
\endaligned$$

Now by Young's inequality, Lemma \ref{temju1} and Lemma \ref{close}, we have
$$\aligned\mathbb{E}\int_0^T|x_\Delta(s)|^{\bar{p}-2}(1+|\bar{x}_\Delta(s)|^r)ds
&\le\mathbb{E}\int_0^T|x_\Delta(s)|^{\bar{p}-2}(1+C(|x_\Delta(s)|^r+|\bar{x}_\Delta(s)-x_\Delta(s)|^r))ds\\&
\le\mathbb{E}\int_0^T(|x_\Delta(s)|^{\bar{p}-2}+C|x_\Delta(s)|^p)ds\\&
\quad+\mathbb{E}\int_0^T|\bar{x}_\Delta(s)-x_\Delta(s)|^pds\le C.\endaligned$$

Since $C$ is independent of $\Delta,$ then the required assertion (\ref{zuida2}) follows. $\square$

\begin{Lemma}\label{jubu1}
Assume that (\ref{local}), (\ref{tiaoj}), (\ref{lianxu}) and (\ref{zengzhang}) hold for $q\le p\le6$. Let $R>|x_0|$ be a positive number and $\Delta_0\le\Delta^*$ be sufficiently small such that $h(\Delta_0)>R$. Let $\theta_{\Delta,R}$ and $e_\Delta(t)$ be the same as before. Then for any $\Delta\in(0,\Delta_0)$ there exists $C$ (independent of $\Delta$ and $R$) such that
$$\mathbb{E}(\sup_{0\le u\le T}|e_\Delta(u\wedge\theta_{\Delta,R})|^{q})<CL^{2q}_{h(\Delta)}\Delta^\frac{q}{2},\ \forall T>0.$$
\end{Lemma}

\textbf{Proof}\ Denote $\theta=\theta_{\Delta,R}$ for simplicity. By It\^o formula, condition (\ref{lianxu}) and BDG inequality again, we have
$$\aligned\mathbb{E}(\sup_{0\le u\le t}|e_\Delta(u\wedge\theta)|^q)&\le q\mathbb{E}\sup_{0\le u\le t}\int_0^{u\wedge\theta}|e_\Delta(s)|^{q-2}
[\langle e_\Delta(s),f(x(s))-f_\Delta(\bar{x}_\Delta(s))\rangle\\&
\qquad+\frac{q-1}{2}|g(x(s))-g_\Delta(\bar{x}_\Delta(s))|^2]ds\\&
\quad+q\mathbb{E}\sup_{0\le u\le t}\left|\int_0^{u\wedge\theta}|e_\Delta(s)|^{q-2}
\langle e_\Delta(s),(g(x(s))-g_\Delta(\bar{x}_\Delta(s)))dB(s)\rangle\right|\\&
\le qH\mathbb{E}\int_0^{t\wedge\theta}|e_\Delta(s)|^{q-2}|x(s)-\bar{x}_\Delta(s)|^2ds\\&
\quad+q\mathbb{E}\int_0^{t\wedge\theta}|e_\Delta(s)|^{q-2}|x_\Delta(s)-\bar{x}_\Delta(s)|
|f(x(s))-f_\Delta(\bar{x}_\Delta(s))|ds\\&
\quad+q\mathbb{E}\sup_{0\le u\le t}\left|\int_0^{u\wedge\theta}|e_\Delta(s)|^{q-2}
\langle e_\Delta(s),(g(x(s))-g_\Delta(\bar{x}_\Delta(s)))dB(s)\rangle\right|.\endaligned$$

Since $\forall 0\le s\le t\wedge\theta,$ $|\bar{x}_\Delta(s)|\le R\le h(\Delta),$ then
$$f_\Delta(\bar{x}_\Delta(s))=f(\bar{x}_\Delta(s)),\quad g_\Delta(\bar{x}_\Delta(s))=g(\bar{x}_\Delta(s)),\quad \forall 0\le s\le t\wedge\theta.$$

Now by BDG inequality again, we have
$$\aligned\mathbb{E}(\sup_{0\le u\le t\wedge\theta}|e_\Delta(u)|^q)&\le  qH\mathbb{E}\int_0^{t\wedge\theta}|e_\Delta(s)|^{q-2}|x(s)-\bar{x}_\Delta(s)|^2ds\\&
\quad+q\mathbb{E}\int_0^{t\wedge\theta}|e_\Delta(s)|^{q-2}|x_\Delta(s)-\bar{x}_\Delta(s)|
|f(x(s))-f(\bar{x}_\Delta(s))|ds\\&\quad+4\sqrt{2}q\mathbb{E}\left|\int_0^{t\wedge\theta}|e_\Delta(s)|^{2q-2}
|g(x(s))-g(\bar{x}_\Delta(s))|^2ds\right|^\frac{1}{2}.\endaligned$$

According to Young's inequality and Lemma \ref{close}, it follows that
$$\aligned\mathbb{E}\int_0^{t\wedge\theta}|e_\Delta(s)|^{q-2}|x(s)-\bar{x}_\Delta(s)|^2ds
&\le\mathbb{E}\int_0^{t\wedge\theta}(\frac{q-2}{q}|e_\Delta(s)|^{q}+\frac{2}{q}|x(s)-\bar{x}_\Delta(s)|^q)ds\\&
\le \mathbb{E}\int_0^{t\wedge\theta}|e_\Delta(s)|^{q}ds+C_{q,T}L^{q}_{h(\Delta)}\Delta^\frac{q}{2}.\endaligned$$

Moreover, as in the proof of Lemma \ref{jubu}, we have
$$\aligned&\quad\mathbb{E}\int_0^{t\wedge\theta}|e_\Delta(s)|^{q-2}|x_\Delta(s)-\bar{x}_\Delta(s)|
|f(x(s))-f(\bar{x}_\Delta(s))|ds\\&\le \mathbb{E}\int_0^{t\wedge\theta}|e_\Delta(s)|^{q-2}
|\bar{x}_\Delta(s)-x_\Delta(s)|\cdot4L_{h(\Delta)}|x(s)-\bar{x}_\Delta(s)|ds\\&
\le C_{1,q}\mathbb{E}\int_0^{t\wedge\theta}|e_\Delta(s)|^{q}ds+C_{2,q,T}L^{2q}_{h(\Delta)}\Delta^\frac{q}{2}.\endaligned$$

Then by (\ref{local}) and Lemma \ref{close}, we have
$$\aligned\mathbb{E}(\sup_{0\le u\le t}|e_\Delta(u\wedge\theta)|^q)&\le   C_{1,q}\mathbb{E}\int_0^{t\wedge\theta}|e_\Delta(s)|^{q}ds+C_{2,q,T}L^{2q}_{h(\Delta)}\Delta^\frac{q}{2}+4\sqrt{2}q
\\&\quad\times\mathbb{E}\left|\sup_{0\le u\le t}|e_\Delta(u\wedge\theta)|^q\int_0^{t\wedge\theta}|e_\Delta(s)|^{q-2}
\cdot16L^2_{h(\Delta)}|x(s)-\bar{x}_\Delta(s)|^2ds\right|^\frac{1}{2}\\&
\le C_{1,q}\mathbb{E}\int_0^{t\wedge\theta}|e_\Delta(s)|^{q}ds+C_{2,q,T}L^{2q}_{h(\Delta)}\Delta^\frac{q}{2}+\frac{1}{2}\mathbb{E}\sup_{0\le u\le t}|e_\Delta(u\wedge\theta)|^q\\&
\quad+16q(q-2)\mathbb{E}\int_0^{t\wedge\theta}|e_\Delta(s)|^{q}ds\\&
\quad+32q\cdot4^q L^q_{h(\Delta)}\int_0^{T}
\mathbb{E}|x(s)-\bar{x}_\Delta(s)|^qds\\&
\le C\int_0^{t}\mathbb{E}\sup_{0\le u\le s}|e_\Delta(u\wedge\theta)|^{q}ds+CL^{2q}_{h(\Delta)}\Delta^\frac{q}{2}.\endaligned$$

Finally, the Gronwall's inequality yields the required assertion. $\square$

\begin{Lemma}\label{close1} Assume that (\ref{local}), (\ref{tiaoj}), (\ref{zengzhang}) and (\ref{gzeng}) hold for $2\le r< p\le6.$ If $q\le p+2-r,$ then for any $\Delta\in(0,\Delta^*)$, there exists $C>0$ (independent of $\Delta$) such that
\begin{equation}\mathbb{E}\left(\sup_{0\le t\le T}|x_\Delta(t)-\bar{x}_\Delta(t)|^q\right)\le CL_{h(\Delta)}^q\Delta^{\frac{q}{2}-1}.\end{equation}
\end{Lemma}
\textbf{Proof}\ Notice that
$$\aligned&\quad\mathbb{E}\sup_{0\le t\le T}|x_\Delta(t)-\bar{x}_\Delta(t)|^q\\&=\mathbb{E}\sup_{0\le k\le[\frac{T}{\Delta}]}\left(\sup_{k\Delta\le t<(k+1)\Delta}|x_\Delta(t)-X_k^\Delta|^q\right)\\&
\le C_q\mathbb{E}\sup_{0\le k\le[\frac{T}{\Delta}]}\left(|f_\Delta(X_k^\Delta)|^q\Delta^q+|g_\Delta(X_k^\Delta)|^q\sup_{k\Delta\le t<(k+1)\Delta}|B(t)-B(k\Delta)|^q\right)\\&
\le C_q\left(4^qL^q_{h(\Delta)}\mathbb{E}\sup_{0\le k\le[\frac{T}{\Delta}]}|X_k^\Delta|^q+|f(0)|^q\right)\Delta^q\\&\quad
+C_q\sum_{k=1}^{[\frac{T}{\Delta}]}\left(4^qL^q_{h(\Delta)}\mathbb{E}|X_k^\Delta|^q+|g(0)|^q\right)\mathbb{E}\sup_{k\Delta\le t<(k+1)\Delta}|B(t)-B(k\Delta)|^q.\endaligned$$

Then by BDG inequality again, we have
$$\aligned\quad\mathbb{E}\sup_{0\le t\le T}|x_\Delta(t)-\bar{x}_\Delta(t)|^q
&\le C_q\left(4^qL^q_{h(\Delta)}\mathbb{E}\sup_{0\le k\le[\frac{T}{\Delta}]}|X_k^\Delta|^q+|f(0)|^q\right)\Delta^q\\&\quad
+C_q\sum_{k=1}^{[\frac{T}{\Delta}]}\left(4^qL^q_{h(\Delta)}\sup_{0\le k\le[\frac{T}{\Delta}]}\mathbb{E}|X_k^\Delta|^q+|g(0)|^q\right)\Delta^\frac{q}{2}\\&
\le C(q,T)L^q_{h(\Delta)}\Delta^{\frac{q}{2}-1}\endaligned$$
as required. $\square$

Now let us prove Theorem \ref{conv1}.

\textbf{Proof of Theorem \ref{conv1}}\ Let $\theta_{\Delta,R}$ and $e_\Delta(t)$ be the same as before.
As in the proof of Theorem \ref{conv},  by Young's inequality, we have that for any $\delta>0,$
$$\aligned\mathbb{E}(\sup_{0\le t\le T}|e_\Delta(t)|^q)&\le \mathbb{E}(1_{\{\theta_{\Delta,R}>T\}}\sup_{0\le t\le T}|e_\Delta(t)|^q)+ \frac{q\delta}{p}\mathbb{E}(\sup_{0\le t\le T}|e_\Delta(t)|^p)\\&\quad+\frac{p-q}{p\delta^{q/(p-q)}}P(\theta_{\Delta,R}\le T).\endaligned$$

By Lemma \ref{zuida}, \ref{zuida1},
$$\mathbb{E}(\sup_{0\le t\le T}|e_\Delta(t)|^p)\le C(\mathbb{E}(\sup_{0\le t\le T}|x(t)|^p)+\mathbb{E}(\sup_{0\le t\le T}|x_\Delta(t)|^p))\le C.$$

Then by Lemma \ref{ting} and \ref{ting1}, we have
$$\aligned\mathbb{E}(\sup_{0\le t\le T}|e_\Delta(t)|^q)&\le \mathbb{E}(\sup_{0\le t\le T}|e_\Delta(t\wedge\theta_{\Delta,R})|^q)+ \frac{Cq\delta}{p}+\frac{C(p-q)}{pR^p\delta^{q/(p-q)}}\endaligned$$
holds for any $\Delta\in(0,\Delta^*), \delta>0$ and $R>|x_0|.$ Choosing
$$\delta=L^{2q}_{h(\Delta)}\Delta^\frac{q}{2},\quad R=(L^{2q}_{h(\Delta)}\Delta^\frac{q}{2})^{-\frac{1}{p-q}},$$
then we get
$$\mathbb{E}(\sup_{0\le t\le T}|e_\Delta(t)|^q)\le\mathbb{E}(\sup_{0\le t\le T}|e_\Delta(t\wedge\theta_{\Delta,R})|^q)+CL^{2q}_{h(\Delta)}\Delta^\frac{q}{2}.$$

Since $h(\Delta)\ge(L^{2q}_{h(\Delta)}\Delta^\frac{q}{2})^{-\frac{1}{p-q}}=R,$ then by Lemma \ref{jubu1}, we have
$$\mathbb{E}(\sup_{0\le t\le T}|e_\Delta(t)|^q)\le CL^{2q}_{h(\Delta)}\Delta^\frac{q}{2}.$$

For the second assertion (\ref{shou2}), we can use (\ref{shou1}) we have just obtained and Lemma \ref{close1} to get
$$\aligned\mathbb{E}\sup_{0\le t\le T}|x(t)-\bar{x}_\Delta(t)|^q&\le C_q\left(\mathbb{E}\sup_{0\le t\le T}|e_\Delta(t)|^q+\mathbb{E}\sup_{0\le t\le T}|x_\Delta(t)-\bar{x}_\Delta(t)|^q\right)\\&
\le C(L^{2q}_{h(\Delta)}\Delta^\frac{q}{2}+L^q_{h(\Delta)}\Delta^{\frac{q}{2}-1}).\endaligned$$

By (\ref{tiaoj}), we have
$$\frac{L^{2q}_{h(\Delta)}\Delta^\frac{q}{2}}{L^q_{h(\Delta)}\Delta^{\frac{q}{2}-1}}=L^q_{h(\Delta)}\Delta\to0$$
since $q\le4.$

Then for sufficiently small $\Delta$, we have
$$\mathbb{E}\sup_{0\le t\le T}|x(t)-\bar{x}_\Delta(t)|^q\le CL^{q}_{h(\Delta)}\Delta^{\frac{q}{2}-1}.$$

We then complete the proof. $\square$

\section{Examples}

Now let us present two examples to illustrate our theory.

\textbf{Example 1}
Consider the following 1-d SDE:
$$dx(t)=(ax(t)-e^{3x(t)})dt+e^{x(t)}dB(t),$$
where $a>0$ and $B(t)$ is a scalar Brownian motion. Then neither $f(x)=ax-e^{3x}$ nor $g(x)=e^x$ is polynomial growing (although both are local Lipschitz continuous). However, we can show that conditions (\ref{zengzhang}) and (\ref{lianxu}) holds for $q=4$ and $p=6.$ Indeed, in this case,
$$xf(x)+\frac{p-1}{2}g^2(x)=ax^2-xe^{3x}+\frac{5}{2}e^{2x}.$$
If $x\le0$ then $-xe^{3x}+\frac{5}{2}e^{2x}=e^{2x}(\frac{5}{2}-xe^x)\le \frac{5}{2}+|x|\le3+x^2$. For $x>0,$ there exists sufficiently large $C>0$ (independent of $x$) such that  $-xe^{3x}+e^{2x}\le C$. So for any $x\in\mathbb{R}^1,$ (\ref{zengzhang}) holds for $p=6$ and $K:=a+C+4.$

On the other hand,
$$\aligned*:&=(x-y)(f(x)-f(y))+\frac{q-1}{2}(g(x)-g(y))^2\\&
=a(x-y)^2-(x-y)(e^{3x}-e^{3y})+\frac{3}{2}(e^{x}-e^{y})^2.\endaligned$$

Choose $R_0>0$ (independent of $x$) sufficiently large.

Without loss of generality, suppose $x>y.$ There are three cases:

\textbf{Case 1}. For any $y<x\le R_0$, according to mean value theorem, we have

$$\aligned *&=a(x-y)^2-3e^{3\theta_1}(x-y)^2+\frac{3}{2}e^{2\theta_2}(x-y)^2\\&\le
(a+\frac{3}{2}e^{2R_0})(x-y)^2,\endaligned$$
where $\theta_i$ lies between $x$ and $y$, $i=1,2.$

\textbf{Case 2}. For $x>y> R_0$, we have
$$\aligned*&=(a+\frac{3}{2}e^{2\theta_2}-3e^{3\theta_1})(x-y)^2\\&
\le(a+\frac{3}{2}e^{2x}-3e^{3y})(x-y)^2.\endaligned$$
If $x-y\le\frac{y}{2}$, then $2x\le3y$ and therefore $*\le a(x-y)^2.$ If $x-y>\frac{y}{2},$ we have
$$\aligned*&=a(x-y)^2-e^{3y}(e^{3(x-y)}-1)(x-y)+\frac{3}{2}e^{2y}(e^{2(x-y)}-1)^2\\&
\le a(x-y)^2+e^{3y}\left(-(x-y)e^{3(x-y)}+(x-y)+\frac{3}{2}(e^{2(x-y)}-2e^{x-y}+1)\right).\endaligned$$
Since in this case $x-y>\frac{y}{2}\ge\frac{R_0}{2}$ (sufficiently large), then
$$-(x-y)e^{3(x-y)}+(x-y)+\frac{3}{2}(e^{2(x-y)}-2e^{x-y}+1)<0.$$

Thus, $*\le a(x-y)^2.$ We have shown that for any $x,y> R_0$, condition (\ref{lianxu}) holds.

\textbf{Case 3}. For $y<R_0<x,$ if $x-y\le R_0,$
then we have
$$*\le a(x-y)^2+\frac{3}{2}e^{2y}(e^{x-y}-1)^2\le a(x-y)^2+\frac{3}{2}e^{2R_0}e^{2\theta}(x-y)^2,$$
where $0\le\theta\le x-y\le R_0.$ Thus, $*\le (a+\frac{3}{2}e^{4R_0})(x-y)^2.$ If $x-y>R_0,$ similar to case 2, we have $*\le a(x-y)^2.$

We have shown that condition (\ref{lianxu}) holds for $H=a+\frac{3}{2}e^{4R_0}$ for any $x,y.$

Moreover, in this case, local Lipschitz condition (\ref{local}) holds for local Lipschitz constant $L_R=3e^{3R}$. Then for any $0<\varepsilon<1,$ we can choose $l(x)=\frac{1}{3^4x^{1-\varepsilon}e^{12x}}$ for $0<x$. It is clear that $l$ is a strict decreasing function in the interval $(0,\infty)$. Let $h$ be the inverse function of $l.$ Then $h$ is also a strict decreasing function in the interval $(0,\Delta^*)$ and $h(\Delta)\to\infty$ as $\Delta\to0$. Now $L^4_{h(\Delta)}\Delta=L^4_Rl(R)=\frac{1}{R^{1-\varepsilon}},$ where $R:=h(\Delta)$. Therefore, $L^4_{h(\Delta)}\Delta=\frac{1}{h(\Delta)^{1-\varepsilon}}\to0$ as $\Delta\to0$. And
$$(L^{2q}_{h(\Delta)}\Delta^{\frac{q}{2}})^{-\frac{1}{p-q}}=(L^{4}_{h(\Delta)}\Delta)^{-1}
 =h(\Delta)^{1-\varepsilon}\le h(\Delta)$$
 for $\Delta$ small enough. Then by Theorem \ref{conv}, for any $T>0$ and sufficient small $\Delta$, we have
\begin{equation}\mathbb{E}|x(T)-x_\Delta(T)|^4\le C(T,2)L^{4}_{h(\Delta)}\Delta=C(T,2)h(\Delta)^{\varepsilon-1}\end{equation}
and
\begin{equation}\mathbb{E}|x(T)-\bar{x}_\Delta(T)|^4\le C(T,2)L^{4}_{h(\Delta)}\Delta=C(T,2)h(\Delta)^{\varepsilon-1}.\end{equation}

Notice that since polynomial growth condition for $f$ fails, then the strong convergence result Theorem 3.4 in \cite{mao1} can not be applied here. However, for the continuous-time MTEM methods (\ref{num1}) and (\ref{num2}), the strong convergence still holds for the given SDE.

\textbf{Example 2}
Consider the scalar SDE
$$dx(t)=(x(t)-x^3(t))dt+|x(t)|^\frac{3}{2}dB(t),$$
where $B(t)$ is a scalar Brownian motion as usual.

In \cite{mao1}, the author only showed that the truncated Euler-Maruyama method strongly converges to the exact solution of the above equation at any fixed time $T$. However, we can prove that the continuous version of modified TEM method (\ref{num2}) strongly converges to the exact solution of the above equation in time interval $[0,T].$

By \cite{mao1}, we know that (\ref{zengzhang}) and (\ref{lianxu}) hold for any $p>3$ and $q=2$. So (\ref{zengzhang}) holds for $p=6$. On the other hand, if we choose $q=4,$ then
$$\aligned &\quad(x-y)(f(x)-f(y))+\frac{q-1}{2}(g(x)-g(y))^2\\&
\le(x-y)^2\left(1-(x^2+y^2+xy)+\frac{9}{4}(|x|^\frac{1}{2}+|y|^\frac{1}{2})\right)\\&
\le (x-y)^2\left(1-\frac{x^2+y^2}{2}+\frac{9}{4}(|x|^\frac{1}{2}+|y|^\frac{1}{2})\right)\\&
\le (x-y)^2\left(\frac{13}{4}-\frac{x^2+y^2}{2}+\frac{9}{8}(|x|+|y|)\right)\\&
\le H(x-y)^2,\endaligned$$
where $H=\frac{13}{4}+\frac{81}{64}.$ That is, (\ref{lianxu}) still holds for $q=4.$

Moreover $|g(x)|^2=|x|^3\le 2(1+|x|^3).$ That is, (\ref{gzeng}) also holds for $r=3$. On the other hand, for any $|x|,|y|\le R,$
$$|x-y-x^3+y^3|\vee||x|^\frac{3}{2}-|y|^\frac{3}{2}|\le ((3R^2+1)\vee\frac{3}{2}\sqrt{R})|x-y|\le (3R^2+1)|x-y|,$$
i.e. $f$ and $g$ are local Lipschitz continuous with local Lipschitz constant $L_R=3R^2+1.$ For $\varepsilon>0$ small enough, choose $$h(x)=\sqrt{\frac{x^\frac{-\varepsilon}{4}-1}{3}},\quad x<1.$$
Then we have $h(\Delta)\to \infty$ and $L^4_{h(\Delta)}\Delta=\Delta^{1-\varepsilon}\to0$ as $\Delta\to0.$ That is, (\ref{tiaoj}) holds for such defined $h$. If we take $\frac{8}{9}<\varepsilon<1,$ then
 $$(L^{2q}_{h(\Delta)}\Delta^{\frac{q}{2}})^{-\frac{1}{p-q}}=(L^{4}_{h(\Delta)}\Delta)^{-1}
 =\Delta^{\varepsilon-1}\le \sqrt{\frac{\Delta^\frac{-\varepsilon}{4}-1}{3}}=h(\Delta),$$
 i.e. (\ref{tj}) holds for small $\Delta$. So by Theorem \ref{conv1}, we have
 $$\mathbb{E}\sup_{0\le t\le T}|x(t)-x_\Delta(t)|^4\le C\Delta^{2(1-\varepsilon)},$$
 and
 $$\mathbb{E}\sup_{0\le t\le T}|x(t)-\bar{x}_\Delta(t)|^4\le CL^q_{h(\Delta)}\Delta^{\frac{q}{2}-1}=C\Delta^{1-\frac{\varepsilon}{2}}.$$
 However, Theorem 4.6 in \cite{mao1} can not be applied here since Assumption 4.1 there does not hold in this case.

 \begin{Remark}
 The above two examples might not be the optimal cases, but they indicate that the modified TEM method needs less conditions than the truncated Euler-Maruyama method introduced by Mao in \cite{mao}.
 \end{Remark}
 
 \section{Conclusions}
 
 We have investigated the strong convergence of so called MTEM method for nonlinear SDE $dX(t)
=f(X(t))dt+g(X(t))dB_t$ in this paper. Strong convergence results are considered for two versions of continuous-time MTEM methods, i.e., continuous-time step-process MTEM method $x_\Delta(t)$ and continuous-time continuous-sample MTEM method $\bar{x}_\Delta(t)$. Strong convergence rates are obtained for both at fixed time $T$ and over a time interval $[0,T]$. For the former, we do not need the polynomial growth condition for $f$, and for the later, we do not need the global Lipschitz condition for $g.$ Therefore, less conditions are needed to ensure the strong convergence for the MTEM method than the truncated EM method.

\end{document}